\documentstyle[amssymb,12pt]{article}

\newtheorem{th}{Theorem}[section]
\newtheorem{lem}[th]{Lemma}
\newtheorem{prop}[th]{Proposition}

\newtheorem{defn}[th]{Definition}
\newenvironment{defn-new}{\begin{defn} \em}{\end{defn}}
\newtheorem{rem}[th]{Remark}
\newenvironment{rem-new}{\begin{rem} \em}{\end{rem}}
\newtheorem{ex}[th]{Example}
\newenvironment{ex-new}{\begin{ex} \em}{\end{ex}}
\newtheorem{exer}[th]{Exercise}
\newenvironment{exer-new}{\begin{exer} \em}{\end{exer}}
\newtheorem{agr}[th]{Agreement}
\newenvironment{agr-new}{\begin{agr} \em}{\end{agr}}
\newtheorem{pbm}[th]{Problem}
\newenvironment{pbm-new}{\begin{pbm} \em}{\end{pbm}}

\makeatletter \@addtoreset{equation}{section} \makeatother

\begin{document}

\begin{center}
%{\Large {\bf ON $(\varepsilon )$-PARA SASAKIAN $3$-MANIFOLDS}} \bigskip \bigskip

{\Large {\bf On $(\varepsilon )$-para Sasakian $3$-manifolds}} \bigskip
\bigskip

Selcen Y\"{u}ksel Perkta\c{s}, Erol K\i l\i \c{c}, Mukut Mani Tripathi and
Sad\i k Kele\c{s}\bigskip \bigskip\bigskip
\end{center}

\noindent {\bf Abstract.} In this paper we study the $3$-dimensional $\left(
\varepsilon \right) $-para Sasakian manifolds. We obtain an necessary and
sufficient condition for an $\left( \varepsilon \right) $-para Sasakian $3$%
-manifold to be an indefinite space form. We show that a
Ricci-semi-symmetric $\left( \varepsilon \right) $-para Sasakian $3$%
-manifold is an indefinite space form. We investigate the necessary and
sufficient condition for an $\left( \varepsilon \right) $-para Sasakian $3$%
-manifold to be locally $\varphi $-symmetric. It is proved that in an $%
\left( \varepsilon \right) $-para Sasakian $3$-manifold with $\eta $%
-parallel Ricci tensor the scalar curvature is constant. It is also shown
that every $\left( \varepsilon \right) $-para Sasakian $3$-manifolds is
pseudosymmetric in the sense of R. Deszcz. \medskip

\noindent {\bf Mathematics Subject Classification:} 53C25, 53C50. \medskip

\noindent {\bf Keywords and phrases:} $\left( \varepsilon \right) $-para
Sasakian $3$-manifold, Ricci-semi-symmetric space, locally $\varphi $%
-symmetric space, $\eta $-parallel Ricci tensor, pseudosymmetric space.

\section{Introduction\label{sect-intro}}

In 1976, S\={a}to \cite{Sato-76} introduced a structure $(\varphi ,\xi ,\eta
)$ satisfying $\varphi ^{2}=I-\eta \otimes \xi $ and $\eta (\xi )=1$ on a
differentiable manifold, which is now well known as an almost paracontact
structure. The structure is an analogue of the almost contact structure \cite%
{Sasaki-60-Tohoku,Blair-02-book} and is closely related to almost product
structure (in contrast to almost contact structure, which is related to
almost complex structure). An almost contact manifold is always
odd-dimensional but an almost paracontact manifold could be even-dimensional
as well. In 1969, T. Takahashi \cite{Takahashi-69-Tohoku-1} introduced
almost contact manifolds equipped with associated pseudo-Riemannian metrics.
In particular, he studied Sasakian manifolds equipped with an associated
pseudo-Riemannian metric. These indefinite almost contact metric manifolds
and indefinite Sasakian manifolds are also known as $\left( \varepsilon
\right) $-almost contact metric manifolds and $\left( \varepsilon \right) $%
-Sasakian manifolds respectively \cite%
{Bej-Dug-93,Duggal-90-IJMMS,Duggal-Sahin-07}. Also, in 1989, K. Matsumoto
\cite{Mat-89} replaced the structure vector field $\xi $ by $-\,\xi $ in an
almost paracontact manifold and associated a Lorentzian metric with the
resulting structure and called it a Lorentzian almost paracontact manifold.
In a Lorentzian almost paracontact manifold given by Matsumoto, the
semi-Riemannian metric has only index $1$ and the structure vector field $%
\xi $ is always timelike. These circumstances motivated the authors in \cite%
{Tripathi-eps-apcm} to associate a semi-Riemannian metric, not necessarily
Lorentzian, with an almost paracontact structure, and they called this
indefinite almost paracontact metric structure an $\left( \varepsilon
\right) $-almost paracontact structure, where the structure vector field $%
\xi $ is spacelike or timelike according as $\varepsilon =1$ or $\varepsilon
=-1$. \medskip

In \cite{Tripathi-eps-apcm} the authors studied $\left( \varepsilon \right) $%
-almost paracontact manifolds, and in particular, $\left( \varepsilon
\right) $-para Sasakian manifolds. They gave basic definitions, some
examples of $\left( \varepsilon \right) $-almost paracontact manifolds and
introduced the notion of an $\left( \varepsilon \right) $-para Sasakian
structure. The basic properties, some typical identities for curvature
tensor and Ricci tensor of the $\left( \varepsilon \right) $-para Sasakian
manifolds were also studied in \cite{Tripathi-eps-apcm}. The authors in \cite%
{Tripathi-eps-apcm} proved that if a semi-Riemannian manifold is one of
flat, proper recurrent or proper Ricci-recurrent, then it can not admit an $%
\left( \varepsilon \right) $-para Sasakian structure. Also they showed that,
for an $\left( \varepsilon \right) $-para Sasakian manifold, the conditions
of being symmetric, semi-symmetric or of constant sectional curvature are
all identical. \medskip

In this paper we study $3$-dimensional $\left( \varepsilon \right) $-para
Sasakian manifolds. The paper organized as follows. Section~\ref{sect-prel}
is devoted to the some basic definitions and curvature properties of $\left(
\varepsilon \right) $-para Sasakian manifolds. In section~\ref{sect-prel},
we also prove that an $\left( \varepsilon \right) $-para Sasakian manifold
is an indefinite space form if and only if the scalar curvature $r$ of the
manifold is equal to $-6\varepsilon $. In section~\ref{sect-eps-R(X,Y).S=0},
we show that a Ricci-semi-symmetric $\left( \varepsilon \right) $-para
Sasakian $3$-manifold is an indefinite space form. In section~\ref%
{sect-eps-phi-symm}, a necessary and sufficient condition for an $\left(
\varepsilon \right) $-para Sasakian $3$-manifold to be locally $\varphi $%
-symmetric is obtained. Section~\ref{sect-eps-eta parallel} contains some
results on $\left( \varepsilon \right) $-para Sasakian $3$-manifolds with $%
\eta $-parallel Ricci tensor. In last section~\ref{sect-eps-PS-ps}, it is
shown that every $\left( \varepsilon \right) $-para Sasakian $3$-manifolds
is pseudosymmetric in the sense of R. Deszcz.

\section{Preliminaries\label{sect-prel}}

Let $M$ be an $n$-dimensional almost paracontact manifold \cite{Sato-76}
equipped with an almost paracontact structure $(\varphi ,\xi ,\eta )$
consisting of a tensor field $\varphi $ of type $(1,1)$, a vector field $\xi
$ and a $1$-form $\eta $ satisfying
\begin{equation}
\varphi ^{2}=I-\eta \otimes \xi ,  \label{eq-phi-eta-xi}
\end{equation}%
\begin{equation}
\eta (\xi )=1,  \label{eq-eta-xi}
\end{equation}%
\begin{equation}
\varphi \xi =0,  \label{eq-phi-xi}
\end{equation}%
\begin{equation}
\eta \circ \varphi =0.  \label{eq-eta-phi}
\end{equation}

Throughout this paper we assume that $X,Y,Z,U,V,W\in {\frak X}\left(
M\right) $, where ${\frak X}\left( M\right) $ is the Lie algebra of vector
fields in $M$, unless specifically stated otherwise. By a semi-Riemannian
metric \cite{ONeill-83} on a manifold $M$, we understand a non-degenerate
symmetric tensor field $g$ of type $\left( 0,2\right) $. In particular, if
its index is $1$, it becomes a Lorentzian metric \cite{Beem-Ehrlich-81}. Let
$g$ be a semi-Riemannian metric with ${\rm index}(g)=\nu $ in an $n$%
-dimensional almost paracontact manifold \ $M$ such that
\begin{equation}
g\left( \varphi X,\varphi Y\right) =g\left( X,Y\right) -\varepsilon \eta
(X)\eta \left( Y\right) ,  \label{eq-metric-1}
\end{equation}%
where $\varepsilon =\pm 1$. Then $M$ is called an $\left( \varepsilon
\right) $-almost paracontact metric manifold equipped with an $\left(
\varepsilon \right) ${\em -}almost paracontact metric structure $(\varphi
,\xi ,\eta ,g,\varepsilon )$ \cite{Tripathi-eps-apcm}. In particular, if $%
{\rm index}(g)=1$, then an $(\varepsilon )$-almost paracontact metric
manifold will be called a Lorentzian almost paracontact manifold. In
particular, if the metric $g$ is positive definite, then an $(\varepsilon )$%
-almost paracontact metric manifold is the usual almost paracontact
metricmanifold \cite{Sato-76}. \medskip

The equation (\ref{eq-metric-1}) is equivalent to
\begin{equation}
g\left( X,\varphi Y\right) =g\left( \varphi X,Y\right)  \label{eq-metric-2}
\end{equation}%
along with
\begin{equation}
g\left( X,\xi \right) =\varepsilon \eta (X).  \label{eq-metric-3}
\end{equation}%
From (\ref{eq-metric-3}) it follows that
\begin{equation}
g\left( \xi ,\xi \right) =\varepsilon ,  \label{eq-g(xi,xi)}
\end{equation}%
that is, the structure vector field $\xi $ is never lightlike. Defining
\begin{equation}
\Phi \left( X,Y\right) \equiv g\left( X,\varphi Y\right) ,
\label{eq-eps-fundamental form}
\end{equation}%
we note that
\begin{equation}
\Phi \left( X,\xi \right) =0.  \label{eq-eps-Phi(X,xi)}
\end{equation}%
\medskip

Let $(M,\varphi ,\xi ,\eta ,g,\varepsilon )$ be an $(\varepsilon )$-almost
paracontact metric manifold (resp. a Lorentzian almost paracontact
manifold). If $\varepsilon =1$, then $M$ will be said to be a spacelike $%
(\varepsilon )$-almost paracontact metric manifold (resp. a spacelike
Lorentzian almost paracontact manifold). Similarly, if $\varepsilon =-\,1$,
then $M$ will be said to be a timelike $(\varepsilon )$-almost paracontact
metric manifold (resp. a timelike Lorentzian almost paracontact manifold)
\cite{Tripathi-eps-apcm}. Note that a timelike Lorentzian almost paracontact
structure is a Lorentzian almost paracontact structure in the sense of Mihai
and Rosca \cite{Mih-Rosca-92,Mat-Mih-Rosca-95}, which differs in the sign of
the structure vector field of the Lorentzian almost paracontact structure
given by Matsumoto \cite{Mat-89}. \medskip

An $\left( \varepsilon \right) $-almost paracontact metric structure is
called an $\left( \varepsilon \right) ${\em -}para Sasakian structure if
\begin{equation}
\left( \nabla _{X}\varphi \right) Y=-\,g(\varphi X,\varphi Y)\xi
-\varepsilon \eta \left( Y\right) \varphi ^{2}X,  \label{eq-eps-PS-def-1}
\end{equation}%
where $\nabla $ is the Levi-Civita connection with respect to $g$. A
manifold endowed with an $\left( \varepsilon \right) $-para Sasakian
structure is called an $\left( \varepsilon \right) $-para Sasakian manifold
\cite{Tripathi-eps-apcm}. In an $\left( \varepsilon \right) $-para Sasakian
manifold we have \cite{Tripathi-eps-apcm}%
\begin{equation}
\nabla \xi =\varepsilon \varphi ,  \label{eq-s-pcm-def}
\end{equation}%
\begin{equation}
\Phi \left( X,Y\right) =g\left( \varphi X,Y\right) =\varepsilon g\left(
\nabla _{X}\xi ,Y\right) =\left( \nabla _{X}\eta \right) Y.
\label{eq-s-pcm-def1}
\end{equation}

An $\left( \varepsilon \right) $-almost paracontact metric manifold is
called $\eta $-Einstein if its Ricci tensor $S$ satisfies the condition%
\begin{equation}
S(X,Y)=ag(X,Y)+b\eta (X)\eta (Y).  \label{eq-eps-PS-eta-Einstein}
\end{equation}%
The $k$-nullity distribution $N(k)$ of a semi-Riemannian manifold $M$ is
defined by%
\begin{equation}
N(k):p\rightarrow N_{p}(k)=\{Z\in T_{p}M:R(X,Y)Z=k(g(Y,Z)X-g(X,Z)Y),
\label{eq-eps-PS-nullity}
\end{equation}%
for all $X,Y\in {\frak X}\left( M\right) $, where $k$ is some smooth
function (see \cite{Tanno}). If $M$ is an $\eta $-Einstein $\left(
\varepsilon \right) $-para Sasakian manifold and the structure vector field $%
\xi $ belongs to the $k$-nullity distribution $N(k)$ for some smooth
function $k$, then we say that $M$ is an $N(k)$-$\eta $-Einstein $\left(
\varepsilon \right) $-para Sasakian manifold (see \cite{Tripathi-Kim}).
\medskip

In an $\left( \varepsilon \right) $-para Sasakian manifold, the Riemann
curvature tensor $R$ and the Ricci tensor $S$ satisfy the following
equations \cite{Tripathi-eps-apcm}:%
\begin{equation}
R\left( X,Y\right) \xi =\eta \left( X\right) Y-\eta \left( Y\right) X,
\label{eq-eps-PS-R(X,Y)xi}
\end{equation}%
\begin{equation}
R\left( X,Y,Z,\xi \right) =-\,\eta \left( X\right) g\left( Y,Z\right) +\eta
\left( Y\right) g\left( X,Z\right) ,  \label{eq-eps-PS-R(X,Y,Z,xi)}
\end{equation}%
\begin{equation}
\eta \left( R\left( X,Y\right) Z\right) =-\,\varepsilon \eta \left( X\right)
g\left( Y,Z\right) +\varepsilon \eta \left( Y\right) g\left( X,Z\right) ,
\label{eq-eps-PS-eta(R(X,Y),Z)}
\end{equation}%
\begin{equation}
R\left( \xi ,X\right) Y=-\,\varepsilon g\left( X,Y\right) \xi +\eta \left(
Y\right) X,  \label{eq-eps-PS-R(xi,X)Y}
\end{equation}%
\begin{equation}
S(X,\xi )=-(n-1)\eta (X).  \label{eq-eps-PS-S(X,xi)}
\end{equation}

It is known that in a semi-Riemannian $3$-manifold
\begin{eqnarray}
R(X,Y)Z &=&g(Y,Z)QX-g(X,Z)QY+S(Y,Z)X-S(X,Z)Y  \nonumber \\
&&-\frac{r}{2}\left( g(Y,Z)X-g(X,Z)Y\right) ,  \label{eq-eps-PS-R(X,Y)Z-3dim}
\end{eqnarray}%
where $Q$ is the Ricci operator and $r$ is the scalar curvature of the
manifold. If we substitute $Z$ by $\xi $ in (\ref{eq-eps-PS-R(X,Y)Z-3dim})
and use (\ref{eq-eps-PS-R(X,Y)xi}), we get%
\begin{equation}
\varepsilon (\eta (Y)QX-\eta (X)QY)=\left( 1+\frac{\varepsilon r}{2}\right)
(\eta (Y)X-\eta (X)Y).  \label{eq-eps-PS-R(X,Y)Z-3dim-Z=xi}
\end{equation}%
By putting $Y=\xi $ in (\ref{eq-eps-PS-R(X,Y)Z-3dim-Z=xi}) and using (\ref%
{eq-eta-xi}) and (\ref{eq-eps-PS-S(X,xi)}) for $n=3$, we obtain%
\[
QX=\frac{1}{2}\{(r+2\varepsilon )X-(r+6\varepsilon )\eta (X)\xi \},
\]%
that is,%
\begin{equation}
S(X,Y)=\frac{1}{2}\{(r+2\varepsilon )g(X,Y)-\varepsilon (r+6\varepsilon
)\eta (X)\eta (Y)\}.  \label{eq-eps-PS-S(X,Y)}
\end{equation}
By using (\ref{eq-eps-PS-S(X,Y)}) in (\ref%
{eq-eps-PS-R(X,Y)Z-3dim}), we obtain%
\begin{eqnarray}
R(X,Y)Z &=&\left( \frac{r}{2}+2\varepsilon \right)
\{g(Y,Z)X-g(X,Z)Y\}
\nonumber \\
&&-\left( \frac{r}{2}+3\varepsilon \right) \{g(Y,Z)\eta (X)\xi
-g(X,Z)\eta
(Y)\xi  \nonumber \\
&&+\varepsilon \eta (Y)\eta (Z)X-\varepsilon \eta (X)\eta (Z)Y\}.
\label{eq-eps-PS-R(X,Y)Z-3dim-final}
\end{eqnarray}%

If an $\left( \varepsilon \right) $-para Sasakian manifold is a space of
constant curvature then it is an indefinite space form.

\begin{lem}
\label{indefinite-spaace-f} An $\left( \varepsilon \right) $-para Sasakian $%
3 $-manifold is an indefinite space form if and only if the scalar curvature
$r=-6\varepsilon .$
\end{lem}

{\em Proof. }Let a $3$-dimensional $\left( \varepsilon \right) $-para
Sasakian manifold be an indefinite space form. Then
\begin{equation}
R(X,Y)Z=c\{g(Y,Z)X-g(X,Z)Y\},\qquad X,Y,Z\in {\frak X}\left( M\right) ,
\label{eq-eps-PS-proof1}
\end{equation}%
where $c$ is the constant curvature of the manifold. By using the definition
of Ricci curvature and (\ref{eq-eps-PS-proof1}) we have%
\begin{equation}
S(X,Y)=2c\,g(X,Y).  \label{eq-eps-PS-proof2}
\end{equation}%
If we use (\ref{eq-eps-PS-proof2}) in the definition of the scalar curvature
we get
\begin{equation}
r=6\,c.  \label{eq-eps-PS-proof3}
\end{equation}%
From (\ref{eq-eps-PS-proof2}) and (\ref{eq-eps-PS-proof3}) one can easily
see that%
\begin{equation}
S(X,Y)=\frac{r}{3}\,g(X,Y).  \label{eq-eps-PS-proof4}
\end{equation}%
By putting $X=Y=\xi $ in (\ref{eq-eps-PS-S(X,Y)}) and using (\ref%
{eq-eps-PS-proof4}) we obtain%
\[
r=-\,6\,\varepsilon .
\]

Conversely, if $ r=-6\varepsilon $ then from the equation (\ref%
{eq-eps-PS-R(X,Y)Z-3dim-final}) we can easily see that the manifold is an
indefinite space form. This completes the proof. $\blacksquare $

\begin{th}
Every $\left( \varepsilon \right) $-para Sasakian $3$-manifold is an $%
N(-\varepsilon )$-$\eta $-Einstein manifold.
\end{th}

{\em Proof.} The proof follows from (\ref{eq-eps-PS-S(X,Y)}) and (\ref%
{eq-eps-PS-R(X,Y)xi}). $\blacksquare $

\section{Ricci-semi-symmetric ($\protect\varepsilon $)-para Sasakian $3$%
-Manifolds \label{sect-eps-R(X,Y).S=0}}

A semi-Riemannian manifold $M$ is said to be Ricci-semi-symmetric \cite%
{Mir-92} if its Ricci tensor $S$ satisfies the condition
\begin{equation}
R(X,Y)\cdot S=0,\qquad X,Y\in {\frak X}\left( M\right) ,
\label{eps-R(X,Y).S=0-1}
\end{equation}%
where $R(X,Y)$ acts as a derivation on $S$. \medskip

Let $M$ be a Ricci-semi-symmetric $(\varepsilon )$-para Sasakian $3$%
-manifold. \ From (\ref{eps-R(X,Y).S=0-1}) we have%
\begin{equation}
S(R(X,Y)U,V)+S(U,R(X,Y)V)=0.  \label{eps-R(X,Y).S=0-2}
\end{equation}%
If we put $Y=\xi $ and use (\ref{eq-eps-PS-R(xi,X)Y})$,$ then we get%
\begin{equation}
\varepsilon g(X,U)S(\xi ,V)-\eta (U)S(X,V)+\varepsilon g(X,V)S(U,\xi )-\eta
(V)S(U,X)=0.  \label{eps-R(X,Y).S=0-3}
\end{equation}%
By using (\ref{eq-eps-PS-S(X,xi)}) in (\ref{eps-R(X,Y).S=0-3}) we obtain%
\begin{equation}
2\varepsilon g(X,U)\eta (V)+S(X,V)\eta (U)+2\varepsilon g(X,V)\eta
(U)+S(X,U)\eta (V)=0.  \label{eps-R(X,Y).S=0-5}
\end{equation}

Consider that $\{e_{1},e_{2},e_{3}\}$ be an orthonormal basis of the $T_{p}M$%
, $p\in M$. Then by putting $X=U=e_{i}$ in (\ref{eps-R(X,Y).S=0-5}) and
taking the summation for $1\leq i\leq 3,$ we have%
\begin{equation}
S(\xi ,V)+8\varepsilon \eta (V)+r\eta (V)=0.  \label{eps-R(X,Y).S=0-6}
\end{equation}%
Again by using (\ref{eq-eps-PS-S(X,xi)}) in (\ref{eps-R(X,Y).S=0-6}), we get%
\[
(r+6\varepsilon )\eta (V)=0,
\]%
which gives $r=-6\varepsilon .$ This implies, in view of Lemma \ref%
{indefinite-spaace-f}, that the manifold is an indefinite space form.
\medskip

Therefore, we can state the following

\begin{th}
A Ricci-semi-symmetric $(\varepsilon )$-para Sasakian $3$-manifold is an
indefinite space form.
\end{th}

\section{Locally $\protect\varphi $-Symmetric $(\protect\varepsilon )$-para
Sasakian $3$-Manifolds \label{sect-eps-phi-symm}}

Analogous to the notion introduced by Takahashi \cite{Takahashi-77-Tohoku}
for Sasakian manifolds, we give the following definition.

\begin{defn-new}
An $(\varepsilon )$-para Sasakian manifold is said to be locally $\varphi $%
-symmetric if%
\[
\varphi ^{2}(\nabla _{W}R)(X,Y,Z)=0,
\]%
for all vector fields $X,$ $Y,$ $Z$ orthogonal to $\xi .$
\end{defn-new}

Now by taking covariant derivative of (\ref{eq-eps-PS-R(X,Y)Z-3dim-final})
with respect to $W$ and using (\ref{eq-eps-fundamental form}) and (\ref%
{eq-eps-Phi(X,xi)}) we have%
\begin{eqnarray}
(\nabla _{W}R)(X,Y,Z) &=&\frac{1}{2}(\nabla _{W}r)\left\{
g(Y,Z)X-g(X,Z)Y\right.  \nonumber \\
&&\qquad \qquad \left. -g(Y,Z)\eta (X)\xi +g(X,Z)\eta (Y)\xi \right.
\nonumber \\
&&\qquad \qquad \left. -\varepsilon \eta (Y)\eta (Z)X+\varepsilon \eta
(X)\eta (Z)Y\right\}  \nonumber \\
&&+\left( \frac{r}{2}+3\varepsilon \right) \left\{ -g(Y,Z)\left( \Phi
(X,W)\xi +\varepsilon \eta (X)\varphi W\right) \right.  \nonumber \\
&&\qquad \qquad +g(X,Z)\left( \Phi (Y,W)\xi +\varepsilon \eta (Y)\varphi
W\right)  \nonumber \\
&&\qquad \qquad -\varepsilon \left( \Phi (Y,W)\eta (Z)+\Phi (Z,W)\eta
(Y)\right) X  \nonumber \\
&&\qquad \qquad \left. +\varepsilon \left( \Phi (X,W)\eta (Z)+\Phi (Z,W)\eta
(X)\right) Y\right\} .  \label{eps-phi-symm-1}
\end{eqnarray}%
Then by taking $X,$ $Y,$ $Z$ orthogonal to $\xi $ and using (\ref%
{eq-phi-eta-xi}), (\ref{eq-phi-xi}), (\ref{eq-eta-phi}) and (\ref%
{eq-metric-3}), from (\ref{eps-phi-symm-1}) we obtain%
\begin{equation}
\varphi ^{2}(\nabla _{W}R)(X,Y,Z)=\frac{1}{2}(\nabla _{W}r)\left(
g(Y,Z)X-g(X,Z)Y\right) .  \label{eps-phi-symm-2}
\end{equation}

Hence from (\ref{eps-phi-symm-2}) we can state the following theorem:

\begin{th}
\label{th-eps-locally-phi-r constant} A $3$-dimensional $(\varepsilon )$%
-para Sasakian manifold is locally $\varphi $-symmetric if and only if the
scalar curvature $r$ is constant.
\end{th}

If a $3$-dimensional $(\varepsilon )$-para Sasakian manifold is
Ricci-semi-symmetric then we have showed that $r=-6\varepsilon $ that is $r$
is constant. Therefore from (\ref{eps-phi-symm-2}), we have

\begin{th}
A $3$-dimensional Ricci-semi-symmetric $(\varepsilon )$-para Sasakian
manifold is locally $\varphi $-symmetric.
\end{th}

In particular, by taking $Z=\xi $ in (\ref{eps-phi-symm-1}) we have
\begin{eqnarray}
(\nabla _{W}R)(X,Y,\xi ) &=&\left( \frac{\varepsilon r}{2}+3\right) \{-\eta
(Y)\Phi (X,W)\xi +\eta (X)\Phi (Y,W)\xi  \nonumber \\
&&\qquad \qquad \quad -\,\Phi (Y,W)X+\Phi (X,W)Y\}.  \label{eps-phi-symm-3}
\end{eqnarray}%
Applying $\varphi ^{2}$ to the both sides of (\ref{eps-phi-symm-3}) we get%
\begin{equation}
\varphi ^{2}(\nabla _{W}R)(X,Y,\xi )=\left( \frac{\varepsilon r}{2}+3\right)
\{-\Phi (Y,W)\varphi ^{2}X+\Phi (X,W)\varphi ^{2}Y\}.  \label{eps-phi-symm-4}
\end{equation}%
If we take $X,$ $Y,$ $W$ orthogonal to $\xi $ in (\ref{eps-phi-symm-3}) and (%
\ref{eps-phi-symm-4}) we have
\[
\varphi ^{2}(\nabla _{W}R)(X,Y,\xi )=(\nabla _{W}R)(X,Y,\xi ).
\]%
Now we can state the following:

\begin{th}
Let $M$ be an $(\varepsilon )$-para Sasakian $3$-manifold such that
\[
\varphi ^{2}(\nabla _{W}R)(X,Y,\xi )=0
\]%
for all $X,Y,W\in {\frak X}\left( M\right) $, orthogonal to $\xi $. Then $M$
is an indefinite space form.
\end{th}

\section{$\protect\eta $-Parallel $(\protect\varepsilon )$-para Sasakian $3$%
-Manifolds \label{sect-eps-eta parallel}}

Motivated by the definitions of Ricci $\eta $-parallelity for Sasakian
manifolds and {\em LP}-Sasakian manifolds were given by Kon \cite{Kon} and
Shaikh and De \cite{Shaikh-De-00}, respectively, we give the following

\begin{defn-new}
Let $M$ be an ($\varepsilon )$-para Sasakian manifold. If the Ricci tensor $%
S $ satisfies
\begin{equation}
(\nabla _{X}S)(\varphi Y,\varphi Z)=0,\qquad X,Y,W\in {\frak X}\left(
M\right) ,  \label{eps-eta-parallel-1}
\end{equation}%
then the manifold $M$ is said to be $\eta $-parallel.
\end{defn-new}

\begin{prop}
\label{pro-eps-eta-parallel-Ricci-r constant} Let $M$ be an $(\varepsilon )$%
-para Sasakian $3$-manifold with $\eta $-parallel Ricci tensor. Then the
scalar curvature $r$ is constant.
\end{prop}

\noindent {\em Proof. } From (\ref{eq-eps-PS-S(X,Y)}) by using (\ref%
{eq-metric-1}) and (\ref{eq-eta-phi})%
\begin{equation}
S(\varphi X,\varphi Y)=\left( \frac{r}{2}+\varepsilon \right) \left(
g(X,Y)-\varepsilon \eta (X)\eta (Y)\right) .  \label{eps-eta-parallel-2}
\end{equation}%
If we take the covariant derivative of (\ref{eps-eta-parallel-2}) with
respect to $Z$ and (\ref{eq-s-pcm-def1}), we get%
\begin{eqnarray}
(\nabla _{Z}S)(\varphi X,\varphi Y) &=&\frac{1}{2}\left\{ (\nabla
_{Z}r)\left( g(X,Y)-\varepsilon \eta (X)\eta (Y)\right) \right.  \nonumber \\
&&\left. -\varepsilon (r+2\varepsilon )\left( \Phi (X,Z)\eta (Y)+\Phi
(Y,Z)\eta (X)\right) \right\} .  \label{eps-eta-parallel-3}
\end{eqnarray}%
Since $M$ is an $(\varepsilon )$-para Sasakian $3$-manifold with $\eta $%
-parallel Ricci tensor, then from (\ref{eps-eta-parallel-1}) we have%
\begin{equation}
(\nabla _{Z}r)\{g(X,Y)-\varepsilon \eta (X)\eta (Y)\}-\varepsilon
(r+2\varepsilon )\{\Phi (X,Z)\eta (Y)+\Phi (Y,Z)\eta (X)\}=0.
\label{eps-eta-parallel-4}
\end{equation}%
Consider that $\{e_{1},e_{2},e_{3}\}$ be an orthonormal basis of the $T_{p}M$%
, $p\in M$. Then by substituting both $X$ and $Y$ by $e_{i},$ $1\leq i\leq 3$%
, in (\ref{eps-eta-parallel-4}) and then taking summation over $i$ and using
(\ref{eq-eps-Phi(X,xi)}) we obtain%
\[
\nabla _{Z}r=0,\qquad Z\in {\frak X}\left( M\right) .
\]%
This completes the proof. $\blacksquare $ \medskip

In view of Theorem \ref{th-eps-locally-phi-r constant} and Proposition \ref%
{pro-eps-eta-parallel-Ricci-r constant} we have the following:

\begin{th}
An $(\varepsilon )$-para Sasakian $3$-manifold with $\eta $-parallel Ricci
tensor is locally $\varphi $-symmetric.
\end{th}

\begin{rem-new}
An $(\varepsilon)$-para Sasakian manifold is called Lorentzian para
Sasakian manifold if $\varepsilon=-1$ and index$(g)=1$. Therefore,
some results we obtained in the previous three sections can be
considered as a generalization of the some results obtained by the
authors in \cite{Shaikh-De-00}.
\end{rem-new}

\begin{rem-new}
In an $(\varepsilon )$-almost para contact $3$-manifold, we observe that $%
{\rm trace}\,\varphi =0$. Therefore the assumption ${\rm
trace}\,\varphi \neq 0$ in \cite{Shaikh-De-00} may not help in
proving several results and some proofs in these papers must be
changed if the results are true anymore.
\end{rem-new}

\section{Pseudosymmetric $(\protect\varepsilon )$-para Sasakian $3$%
-Manifolds \label{sect-eps-PS-ps}}

Now, we consider a well known generalization of the concept of an $\eta $%
-Einstein almost paracontact metric manifold in the following

\begin{defn-new}
\cite{Chaki-Maity-00} A non-flat $n$-dimensional Riemannian manifold $\left(
M,g\right) $ is said to be a {\em quasi Einstein manifold} if its Ricci
tensor $S$ satisfies
\begin{equation}
S=ag+b\eta \otimes \eta  \label{qE-defn-S}
\end{equation}%
or equivalently, its Ricci operator $Q$ satisfies
\begin{equation}
Q=aI+b\eta \otimes \xi  \label{qE-defn}
\end{equation}%
for some smooth functions $a$ and $b$, where $\eta $ is a nonzero $1$-form
such that
\begin{equation}
g\left( X,\xi \right) =\eta \left( X\right) ,\qquad g\left( \xi ,\xi \right)
=\eta \left( \xi \right) =1  \label{eta-xi}
\end{equation}%
for the associated vector field $\xi $. The $1$-form $\eta $ is called the
associated $1$-form and the unit vector field $\xi $ is called the generator
of the quasi Einstein manifold.
\end{defn-new}

B. Y. Chen and K. Yano \cite{Chen-Yano-72} defined a Riemannian manifold $%
\left( M,g\right) $ to be of {\em quasi-constant curvature} if it is
conformally flat manifold and its Riemann-Christoffel curvature tensor $R$
of type $(0,4)$ satisfies the condition
\begin{eqnarray}
R(X,Y,Z,W) &=&a\left\{ g(Y,Z)g(X,W)-g(X,Z)g(Y,W)\right\}  \nonumber \\
&&+\,b\left\{ g(X,W)T(Y)T(Z)-g(X,Z)T(Y)T(W)\right.  \nonumber \\
&&\qquad \left. +\,g(Y,Z)T(X)T(W)-g(Y,W)T(X)T(Z)\right\}
\label{eq-quasi-const-curv}
\end{eqnarray}%
for all $X,Y,Z,W\in {\frak X}\left( M\right) $, where $a$, $b$ are some
smooth functions and $T$ is a non-zero $1$-form defined by
\[
g(X,\rho )=T(X),\qquad X\in {\frak X}\left( M\right)
\]%
for a unit vector field $\rho $. On the other hand, Gh. Vr\u{a}nceanu \cite%
{Vranceanu-68-book} defined a Riemannian manifold $\left( M,g\right) $ to be
of {\em almost constant curvature} if \ $M$ satisfies (\ref%
{eq-quasi-const-curv}). Later on, it was pointed out by A. L. Mocanu \cite%
{Mocanu-87} that the manifold introduced by Chen and Yano and the manifold
introduced by Vr\u{a}nceanu were identical, as it can be verified that if
the curvature tensor $R$ is of the form (\ref{eq-quasi-const-curv}), then
the manifold is conformally flat. Thus, a Riemannian manifold is said to be
of {\em quasi-constant curvature} if the curvature tensor $R$ satisfies (\ref%
{eq-quasi-const-curv}). If $b=0$, then the manifold reduces to a manifold of
constant curvature.

\begin{ex-new}
A manifold of quasi-constant curvature is a quasi Einstein manifold \cite[%
Example 1]{De-Ghosh-04-PMH}. Conversely, a conformally flat quasi Einstein
manifold of dimension $n$ ($n>3$) is a manifold of quasi-constant curvature
\cite[Theorem 4]{De-Ghosh-05-Debrecen}.
\end{ex-new}

Let $\left( M,g\right) $ be a semi-Riemannian manifold with its Levi-Civita
connection $\nabla $. A tensor field $F$ of type $\left( 1,3\right) $ is
known to be {\em curvature-like} provided that $F$ satisfies the symmetric
properties of the curvature tensor $R$. For example, the tensor $R_{g}$
given by
\begin{equation}
R_{g}\left( X,Y\right) Z\equiv \left( X\wedge _{g}Y\right)
Z=g(Y,Z)X-g(X,Z)Y,\qquad X,Y\in {\frak X}\left( M\right) ,  \label{eq-R0}
\end{equation}%
is a trivial example of a curvature like tensor. Sometimes, the symbol $%
R_{g} $ seems to be much more convenient than the symbol $\left( X\wedge
_{g}Y\right) Z$. For example, a semi-Riemannian manifold $\left( M,g\right) $
is of constant curvature $c$ if and only if $R=cR_{g}$. \medskip

It is well known that every curvature-like tensor field $F$ acts on the
algebra ${\frak T}_{s}^{1}\left( M\right) $ of all tensor fields on $M$ of
type $\left( 1,s\right) $ as a derivation \cite[p. 44]{ONeill-83}:
\begin{eqnarray*}
\left( F\cdot P\right) \left( X_{1},\ldots ,X_{s};Y,X\right)
&=&F(X,Y)\left\{ P\left( X_{1},\ldots ,X_{s}\right) \right\} \\
&&\qquad -\sum_{j=1}^{s}P(X_{1},\ldots ,F\left( X,Y\right) X_{j},\ldots
,X_{s})
\end{eqnarray*}%
for all $X_{1},\ldots ,X_{s}\in {\frak X}\left( M\right) $, $P\in {\frak T}%
_{s}^{1}\left( M\right) $. The derivative $F\cdot P$ of $P$ by $F$ is a
tensor field of type $\left( 1,s+2\right) $. A semi-Riemannian manifold $%
\left( M,g\right) $ is said to be {\em semi-symmetric} if $R\cdot R=0$.
Obviously, locally symmetric spaces ($\nabla R=0$) are semi-symmetric. More
generally, a semi-Riemannian manifold $\left( M,g\right) $ is said to be
{\em pseudo-symmetric} (in the sense of R. Deszcz) \cite{Desz-87} if $R\cdot
R$ and $R_{g}\cdot R\/$ in $M$ are linearly dependent, that is, if there
exists a real valued smooth function $L:M\rightarrow {\Bbb R}$ such that
\[
R\cdot R=L\,R_{g}\cdot R
\]%
is true on the set
\[
U=\left\{ x\in M:R\neq \frac{r}{n(n-1)}R_{g}\;{\rm at}\;x\right\} .
\]%
A pseudo-symmetric space is said to be {\em proper} if it is not
semi-symmetric. For details we refer to \cite{Bo-Kow-Van-book,
Belkhelfa-01-Thesis}. \medskip

In the literature, there is also another notion of pseudo-symmetry. A
semi-Riemannian manifold $\left( M,g\right) $ is said to be {\em %
pseudo-symmetric} in the sense of Chaki \cite{Chaki-87-Iasi} if
\begin{eqnarray*}
(\nabla R)(X_{1},X_{2},X_{3},X_{4};X) &=&2\omega
(X)R(X_{1},X_{2},X_{3},X_{4})+\omega (X_{1})R((X,X_{2},X_{3},X_{4}) \\
&&\quad +\,\omega (X_{2})R((X_{1},X,X_{3},X_{4})+\omega
(X_{3})R((X_{1},X_{2},X,X_{4}) \\
&&\quad +\,\omega (X_{4})R((X_{1},X_{2},X_{3},X)
\end{eqnarray*}%
for all $X_{1},X_{2},X_{3},X_{4};X\in {\frak X}\left( M\right) $, where $%
\omega $ is a $1$-form on $(M,g)$. Of course, both the definitions of
pseudo-symmetry for a semi-Riemannian manifold are not equivalent. For
example, in contact geometry, every Sasakian space form is pseudo-symmetric
in the sense of Deszcz \cite[Theorem 2.3]{BDV-05}, but a Sasakian manifold
cannot be pseudo-symmetric in the sense of Chaki \cite[Theorem~1]%
{Tarafdar-91-PMH}. We assume the pseudo-symmetry always in the sense of
Deszcz, unless specifically stated otherwise. \medskip

For Riemannian $3$-manifolds, the following characterization of
pseudosymmetry is known (cf. \cite{Kow-Seki-96,Cho-Inoguchi-Lee-09}).

\begin{prop}
A $3$-dimensional Riemannian manifold $\left( M,g\right) $ is
pseudo-symmetric if and only if it is quasi-Einstein, that is, if and only
if there exists a $1$-form $\eta $ such that the Ricci tensor field $S$
satisfies $S=ag+b\eta \otimes \eta $ for some smooth functions $a$ and $b$.
\end{prop}

In view of the above Proposition, we can state the following:

\begin{th}
Every $3$-dimensional $\eta $-Einstein $\left( \varepsilon \right) $-almost
paracontact metric manifold is always pseudo-symmetric. In particular, each $%
3$-dimensional $\left( \varepsilon \right) $-para Sasakian manifold is
pseudo-symmetric.
\end{th}

\begin{pbm-new}
It would be interesting to know whether an $\left( \varepsilon \right) $%
-almost para Sasakian manifold is pseudo-symmetric in the sense of Chaki or
not.
\end{pbm-new}

\noindent Selcen Y\"{u}ksel Perkta\c{s}

\noindent Department of Mathematics, Faculty of Arts and Sciences, \.{I}n%
\"{o}n\"{u} University

\noindent 44280 Malatya, Turkey

\noindent Email: selcenyuksel@inonu.edu.tr

\smallskip

\noindent Erol K\i l\i \c{c}

\noindent Department of Mathematics, Faculty of Arts and Sciences, \.{I}n%
\"{o}n\"{u} University

\noindent 44280 Malatya, Turkey

\noindent Email: ekilic@inonu.edu.tr

\smallskip

\noindent Mukut Mani Tripathi

\noindent Department of Mathematics, Banaras Hindu University

\noindent Varanasi 221 005, India.

\noindent Email: mmtripathi66@yahoo.com

\smallskip

\noindent Sad\i k Kele\c{s}

\noindent Department of Mathematics, Faculty of Arts and Sciences, \.{I}n%
\"{o}n\"{u} University

\noindent 44280 Malatya, Turkey

\noindent Email: skeles@inonu.edu.tr

\end{document}